\documentclass[11pt,reqno]{amsart}

\usepackage[T1]{fontenc}
\usepackage[utf8]{inputenc}
\usepackage{amsmath,amssymb,amsthm}
\usepackage[letterpaper,margin=0.75in]{geometry}
\usepackage{booktabs}
\usepackage[hidelinks]{hyperref}

\DeclareMathOperator{\vol}{vol}
\DeclareMathOperator{\conv}{conv}
\DeclareMathOperator{\area}{area}
\DeclareMathOperator{\spn}{span}

\newcommand{\repobase}{https://github.com/csabamedgyeswork/inscribed-polytopes}
\newcommand{\repocommit}{maxvol-d4-n9-v1}
\newcommand{\repohref}{\href{\repobase}{\texttt{\repobase}}}
\newcommand{\artpath}[1]{\href{\repobase/tree/\repocommit/#1}{\nolinkurl{#1}}}

\newcounter{artifact}
\newenvironment{artifactlist}
  {\sloppy 
   \begin{list}{[R\arabic{artifact}]}
     {\usecounter{artifact}%
      \setlength{\labelwidth}{3.1em}\setlength{\leftmargin}{3.7em}%
      \setlength{\labelsep}{0.6em}\setlength{\itemsep}{0.15ex}%
      \setlength{\parsep}{0pt}}}
  {\end{list}}
\newcommand{\artifact}[3]{\item\label{art:#1}\artpath{#2}, #3}
\newcommand{\artref}[1]{\hyperref[art:#1]{[R\ref*{art:#1}]}}

\newtheorem{theorem}{Theorem}[section]
\newtheorem{proposition}[theorem]{Proposition}
\newtheorem{lemma}[theorem]{Lemma}
\newtheorem{corollary}[theorem]{Corollary}

\theoremstyle{definition}
\newtheorem{definition}[theorem]{Definition}
\newtheorem{remark}[theorem]{Remark}
\newtheorem*{notation}{Notation}

\newcommand{\citedname}{}
\newtheorem*{citedstatement}{\citedname}
\newenvironment{cited}[2]
  {\renewcommand{\citedname}{#1}\begin{citedstatement}[#2]}
  {\end{citedstatement}}

\title[The inscribed $4$-polytope of maximal volume with nine vertices]
      {The inscribed $4$-polytope of maximal volume\\ with nine vertices}
\author{Csaba Medgyes}
\email{csaba.medgyes.work@gmail.com}

\begin{document}

\begin{abstract}
  We prove that a $4$-polytope with nine vertices inscribed in the unit sphere
  $\mathbb{S}^{3}$ has volume at most $\frac{5}{24}\sqrt{10 + 2\sqrt{5}}$, and
  that, up to orthogonal transformations, equality holds only for the convex
  hull of a square and a regular pentagon inscribed in $\mathbb{S}^{3}$ in
  orthogonal $2$-planes through the origin. For the combinatorial type of this
  polytope, whose boundary complex is the join of a $4$-cycle and a $5$-cycle,
  we write the signed volume as the pairing,
  through the Hodge star, of the area vectors of the two cycles, and bound each
  area vector by finite Fourier analysis; the bound holds for every choice of
  nine points on the sphere. Since an optimal polytope is
  simplicial, by a result of Horv\'ath and L\'angi, it remains to exclude the
  other simplicial types, which we do by sum-of-squares certificates checked in
  exact arithmetic.
\end{abstract}

\maketitle

\section{Introduction}
\label{sec:intro}

Which convex polytope with $n$ vertices, inscribed in the unit sphere
$\mathbb{S}^{d-1}$ of $\mathbb{R}^{d}$, has the largest volume? Fejes T\'{o}th
raised this question in three dimensions \cite{fejes-toth-1964}, and it is
listed among the open problems of discrete geometry
\cite{brass-moser-pach-2005}. There, besides the classical cases $n = 4, 6,
12$ (see \cite[Section~1.2]{hoehner-ledford-2026}), Berman and Hanes
determined the optimal polyhedra with at most eight vertices
\cite{berman-hanes-1970}, and Hoehner and Ledford recently those with nine
\cite{hoehner-ledford-2026}. For more vertices only numerical results are
known \cite{mutoh-2003}, \cite{hardin-sloane-smith}.

In higher dimensions even less is known. Horv\'{a}th and L\'{a}ngi extended
the necessary condition for optimality of Berman and Hanes to every dimension,
and with it showed that optimal polytopes are simplicial
\cite[Lemma~1]{horvath-langi-2016}. They determined the optimal polytopes with
$n = d + 1$ vertices, the regular simplices
\cite[Corollary~1]{horvath-langi-2016}; with $n = d + 2$, the convex hulls of
two regular simplices of dimensions $\lfloor d/2\rfloor$ and $\lceil
d/2\rceil$ in orthogonal subspaces \cite[Theorem~2]{horvath-langi-2016}; and,
for odd $d$, with $n = d + 3$ \cite[Theorem~3]{horvath-langi-2016}. For even
$d$ their argument for $n = d + 3$ excludes cyclic polytopes, and they asked
whether an optimal polytope can be cyclic
\cite[Problem~2]{horvath-langi-2016}.

In dimension four these results settle $n = 5$ and $n = 6$. In earlier work
we answered the question of Horv\'{a}th and L\'{a}ngi for $d = 4$, which
settles $n = 7$, and showed that the regular cross-polytope is optimal for
$n = 8$ \cite{medgyes-maxvol-d4}. The optima for $n = 6$, $7$ and $8$ share a
structure: each is the convex hull of two regular polygons with $m$ and $l$
vertices, inscribed in $\mathbb{S}^{3}$ in orthogonal $2$-planes through the
origin, with $(m, l) = (3, 3)$, $(3, 4)$ and $(4, 4)$. To our knowledge, no case with $d = 4$
and $n \geq 9$ has been settled.

This paper treats the first open case in dimension four, $n = 9$, where the
pattern suggests a square and a regular pentagon. The method of
\cite{medgyes-maxvol-d4} bounds the volume on each combinatorial type of
simplicial polytope by a sum-of-squares certificate, checked in exact
arithmetic. Two things change with nine vertices. The number of types grows
from $37$ \cite{grunbaum-sreedharan-1967} to $1142$
\cite{altshuler-bokowski-steinberg-1980}. And the volume of the square--pentagon pair is
irrational, while a certificate gives a rational bound, so on the type of the
pair no certificate is sharp.

We prove that the square--pentagon pair is, up to orthogonal transformations,
the unique optimum (Theorem~\ref{thm:main}). Certificates bound the volume on
every other type strictly below that of the pair. On the type of
the pair, whose boundary complex is the join of a $4$-cycle and a $5$-cycle,
we prove a sharp
inequality instead. Its two ingredients are a factorisation of the signed
volume through the Hodge star, and a bound on the area vector of a closed
polygon with unit vertices, proved by expanding the vertex coordinates in the
real discrete Fourier basis, as Schoenberg did for plane polygons
\cite{schoenberg-1950}.

Section~\ref{sec:prelim} collects the definitions and known results we use.
Section~\ref{sec:lemmas} proves three lemmas about the join of two cycles,
among them the factorisation of the signed volume, and
Section~\ref{sec:polygon} the inequality for closed polygons.
Section~\ref{sec:pair} proves the sharp bound within the join type,
Section~\ref{sec:rivals} bounds all other types, and
Section~\ref{sec:main} combines the two.

\section{Preliminaries}
\label{sec:prelim}

\subsection{Inscribed polytopes}

\begin{cited}{Notation}{{\cite[Section~1]{horvath-langi-2016}}}
  Throughout the paper, by a polytope we mean a convex polytope. The vertex set
  of a polytope $P$ is denoted by $V(P)$. We denote the family of
  $d$-dimensional polytopes, with $n$ vertices and inscribed in the unit sphere
  $\mathbb{S}^{d-1}$, by $\mathcal{P}_{d}(n)$. We denote $d$-dimensional
  volume by $\vol_{d}$, and set
  $v_{d}(n) = \max\{\vol_{d}(P) : P \in \mathcal{P}_{d}(n)\}$.
\end{cited}

\begin{remark}[{\cite[Remark~2.1]{medgyes-maxvol-d4}}]
  \label{rem:attained}
  For $n \geq d + 1$, $v_{d}(n)$ is attained by some
  $P \in \mathcal{P}_{d}(n)$.
\end{remark}

\begin{definition}[{\cite[Definition~2.1, Section~2.2 and Proposition~2.16]{ziegler-1995}}]
  A \emph{facet} of a $d$-polytope is a face of dimension $d-1$, and a
  \emph{ridge} is a face of dimension $d-2$. A polytope is \emph{simplicial} if
  every facet is a simplex. Two polytopes are \emph{combinatorially
  equivalent} if their face lattices are isomorphic; an equivalence class is a
  \emph{combinatorial type}.
\end{definition}

Equivalently, two polytopes are combinatorially equivalent if their
facet--vertex incidence matrices differ only by permutations of rows and
columns \cite[Section~15.1.1]{henk-richter-gebert-ziegler-2017}. For
simplicial polytopes, which is the form in which we use the definition, this
means that a bijection between the vertex sets carries the vertex sets of the
facets of one onto those of the other.

\begin{definition}[{\cite[Definition~1]{horvath-langi-2016}}]
  \label{def:propertyz}
  Let $P \in \mathcal{P}_{d}(n)$ be a $d$-polytope with
  $V(P) = \{p_{1}, p_{2}, \ldots, p_{n}\}$. If for each $i$, there is an open
  set $U_{i} \subset \mathbb{S}^{d-1}$ such that $p_{i} \in U_{i}$, and for any
  $q \in U_{i}$, we have
  \begin{equation*}
    \vol_{d}\bigl(\conv((V(P) \setminus \{p_{i}\}) \cup \{q\})\bigr)
      \leq \vol_{d}(P),
  \end{equation*}
  then we say that $P$ satisfies Property Z.
\end{definition}

\begin{remark}[{\cite{berman-hanes-1970}, \cite[Remark~2.4]{medgyes-maxvol-d4}}]
  \label{rem:propertyz}
  Every $P \in \mathcal{P}_{d}(n)$ with $\vol_{d}(P) = v_{d}(n)$ satisfies
  Property Z.
\end{remark}

\begin{lemma}[{\cite[Lemma~1]{horvath-langi-2016}}]
  \label{lem:simplicial}
  Every $P \in \mathcal{P}_{d}(n)$ satisfying Property Z is simplicial.
\end{lemma}

\subsection{Signed volume and certificates}

For $v_{1}, \ldots, v_{d} \in \mathbb{R}^{d}$, $\det(v_{1}, \ldots, v_{d})$
denotes the determinant of the matrix with columns $v_{1}, \ldots, v_{d}$.

\begin{definition}[{\cite[Definition~3.1]{medgyes-maxvol-d4}}]
  \label{def:signed-volume}
  An \emph{orientation} of a finite set is an ordering of its elements up to
  even permutations. Let $\mathcal{C} = (n, \mathcal{F})$, where
  $\mathcal{F}$ is a finite set of ordered $d$-tuples
  $F = (F_{1}, \ldots, F_{d})$ of distinct elements of $\{1, \ldots, n\}$,
  called facets. Deleting $F_{a}$ gives a ridge, with the orientation of
  $(F_{1}, \ldots, \widehat{F_{a}}, \ldots, F_{d})$ if $a$ is even and the
  opposite one if $a$ is odd. The orientation of $\mathcal{C}$ is
  \emph{coherent} if every ridge lies in exactly two facets and these induce
  opposite orientations on it. A point
  $p = (p_{1}, \ldots, p_{n}) \in (\mathbb{R}^{d})^{n}$ \emph{realises}
  $\mathcal{C}$ if $P = \conv\{p_{1}, \ldots, p_{n}\}$ has the $p_{i}$ as its
  vertices and its facets are exactly the simplices
  $\conv\{p_{j} : j \in F\}$, $F \in \mathcal{F}$. The \emph{signed volume}
  of $p$ is
  \begin{equation}
    \label{eq:signed-volume}
    D(p) = \sum_{F \in \mathcal{F}} \det(p_{F_{1}}, \ldots, p_{F_{d}}).
  \end{equation}
\end{definition}

\begin{lemma}[{\cite[Lemma~3.2]{medgyes-maxvol-d4}}]
  \label{lem:volume}
  If the orientation of $\mathcal{C}$ is coherent and $p$ realises
  $\mathcal{C}$, then
  \[
    d!\,\vol_{d}(\conv\{p_{1}, \ldots, p_{n}\}) = |D(p)|.
  \]
\end{lemma}

\begin{definition}[{\cite[Definition~3.4]{medgyes-maxvol-d4}}]
  \label{def:certificate}
  Let $x_{i} = (x_{i,1}, \ldots, x_{i,d})$ be the coordinates of $p_{i}$, and
  $v(x)$ the vector of all monomials of degree at most $2$ in the $x_{i,a}$. A
  \emph{certificate} for $\mathcal{C}$ and $\beta \in \mathbb{Q}$ is a
  symmetric positive definite matrix $Q$ with rational entries together with
  polynomials $h_{1}, \ldots, h_{n}$ with rational coefficients such that
  \begin{equation}
    \label{eq:certificate}
    \beta - D(x) = v(x)^{\mathsf{T}} Q\, v(x)
      + \sum_{i=1}^{n} h_{i}(x)\bigl(|x_{i}|^{2} - 1\bigr)
  \end{equation}
  holds identically as polynomials in the $x_{i,a}$.
\end{definition}

\begin{proposition}[{\cite[Proposition~3.5]{medgyes-maxvol-d4}}]
  \label{prop:certificate}
  Suppose the orientation of $\mathcal{C}$ is coherent and a certificate for
  $\mathcal{C}$ and $\beta$ exists. Then every polytope inscribed in
  $\mathbb{S}^{d-1}$ that realises $\mathcal{C}$ has volume at most
  $\beta/d!$.
\end{proposition}

A certificate is finite rational data, and its hypotheses (the identity
\eqref{eq:certificate}, positive definiteness of $Q$, coherence, and the facet
sets of the type) are checked in exact arithmetic as described in
\cite[Section~3.3]{medgyes-maxvol-d4}.

\subsection{Simplicial \texorpdfstring{$4$}{4}-polytopes with nine vertices}

\begin{theorem}[{\cite{altshuler-bokowski-steinberg-1980},
    \cite[Table~3]{firsching-2017}}]
  \label{thm:count}
  There are exactly $1142$ combinatorial types of simplicial $4$-polytopes with
  $9$ vertices.
\end{theorem}

The count is due to Altshuler, Bokowski and Steinberg; table numbers of
\cite{firsching-2017} refer to its arXiv version.
We use the realisations of all $1142$ types given there (Lemma~\ref{lem:data}).

\subsection{Direct sums}

\begin{cited}{Direct sum}{{\cite[Section~15.1.3]{henk-richter-gebert-ziegler-2017}}}
  Let $P \subseteq \mathbb{R}^{d}$ be a $d$-polytope with $n$ vertices and $m$
  facets, and $P' \subseteq \mathbb{R}^{d'}$ a $d'$-polytope with $n'$
  vertices and $m'$ facets, both with the origin in their interiors. The
  \emph{direct sum}
  $P \oplus P' = \conv(\{(p, 0) : p \in P\} \cup \{(0, p') : p' \in P'\})$ is
  $(d + d')$-dimensional, and has $n + n'$ vertices and $m \cdot m'$ facets.
\end{cited}

We use this for two regular polygons in orthogonal $2$-planes through the
origin, in Lemma~\ref{lem:freesum}. A polygon has as many edges (its facets)
as vertices, so for polygons with $m$ and $l$ vertices the direct sum has
$m + l$ vertices and $ml$ facets. For $m \geq 3$, $R_{m}$ denotes a regular
$m$-gon inscribed in a unit circle; $\area(R_{m}) = \frac{m}{2}\sin\frac{2\pi}{m}$.

\subsection{Bivectors}

Let $e_{1}, \ldots, e_{N}$ be the standard basis of $\mathbb{R}^{N}$ and
$v \cdot w$ the standard inner product. We use the usual identification of
bivectors with skew-symmetric matrices: the space $\Lambda^{2}\mathbb{R}^{N}$
of \emph{bivectors} is the space of real skew-symmetric $N \times N$ matrices,
and for $v, w \in \mathbb{R}^{N}$ the \emph{wedge product} is
\[
  v \wedge w = vw^{\mathsf{T}} - wv^{\mathsf{T}}.
\]
It is bilinear and satisfies $w \wedge v = -\,v \wedge w$, so $v \wedge v = 0$.
A skew-symmetric matrix $X$ is determined by its entries $X_{kl}$, $k < l$,
above the diagonal, and $e_{kl} = e_{k} \wedge e_{l}$ has entry $1$ at
$(k, l)$, $-1$ at $(l, k)$ and $0$ elsewhere. Hence the $e_{kl}$, $k < l$, form
a basis of $\Lambda^{2}\mathbb{R}^{N}$, we write
$X = \sum_{k<l} X_{kl}\, e_{kl}$, and
\[
  (v \wedge w)_{kl} = v_{k} w_{l} - v_{l} w_{k},
\]
which is the signed area of the projection of the parallelogram spanned by
$v$ and $w$ onto the plane of $e_{k}$ and $e_{l}$. Every bivector is a sum of
wedge products, but for $N \geq 4$ not every bivector is a single wedge
product: $v \wedge w$ has rank at most $2$, while $e_{12} + e_{34}$ has
rank $4$.

We use the inner product
$\langle X, Y\rangle = \frac12\operatorname{tr}(X^{\mathsf{T}}Y)
= \sum_{k<l} X_{kl}Y_{kl}$ on $\Lambda^{2}\mathbb{R}^{N}$, and
$|X| = \langle X, X\rangle^{1/2}$; the factor $\frac12$ accounts for each pair
$k < l$ appearing twice in $\operatorname{tr}(X^{\mathsf{T}}Y) = \sum_{i,j}
X_{ij}Y_{ij}$. Expanding the trace with
$\operatorname{tr}\bigl((pq^{\mathsf{T}})^{\mathsf{T}} rs^{\mathsf{T}}\bigr)
= (p \cdot r)(q \cdot s)$ gives
\[
  \langle a \wedge b, c \wedge d\rangle = (a \cdot c)(b \cdot d) - (a \cdot d)(b \cdot c).
\]
By this formula, $u_{k} \wedge u_{l}$, $k < l$, are orthonormal for every
orthonormal basis $u_{1}, \ldots, u_{N}$ of $\mathbb{R}^{N}$. Taking $c = a$
and $d = b$ in the formula gives the following.

\begin{lemma}[Lagrange's identity]
  \label{lem:lagrange}
  For $a, b \in \mathbb{R}^{N}$,
  $|a \wedge b|^{2} = |a|^{2}|b|^{2} - (a \cdot b)^{2}$.
\end{lemma}

We use Lemma~\ref{lem:lagrange} in the form $|a \wedge b| \leq |a|\,|b| \leq
\frac12(|a|^{2} + |b|^{2})$, where both inequalities are equalities if and only
if $a \cdot b = 0$ and $|a| = |b|$.

Let $\Pi \subseteq \mathbb{R}^{N}$ be a $2$-plane through the origin and
$\Pi^{\perp}$ its orthogonal complement. For an orthonormal basis $u_{1},
u_{2}$ of $\Pi$, the \emph{unit area element} of $\Pi$ is
$\omega_{\Pi} = u_{1} \wedge u_{2}$; for instance, $e_{kl}$ is that of the
plane spanned by $e_{k}$ and $e_{l}$. By bilinearity and antisymmetry,
changing the orthonormal basis multiplies $u_{1} \wedge u_{2}$ by the
determinant of the change of basis, which is $\pm 1$; statements about
$\omega_{\Pi}$ below hold for either sign. By Lemma~\ref{lem:lagrange}, $|\omega_{\Pi}| = 1$.

\subsection{The Hodge star on \texorpdfstring{$\Lambda^{2}\mathbb{R}^{4}$}{Lambda2R4}}

We use the Hodge star only for $N = 4$. For $k < l$ in $\{1, 2, 3, 4\}$ let
$\{k', l'\}$, $k' < l'$, be the complement of $\{k, l\}$, so that the plane
spanned by $e_{k'}, e_{l'}$ is the orthogonal complement of the plane spanned
by $e_{k}, e_{l}$. Let $\varepsilon_{kl} = \pm 1$ be the sign of the
permutation $(k, l, k', l')$ of $(1, 2, 3, 4)$. The \emph{Hodge star} is the
linear map $\star$ of $\Lambda^{2}\mathbb{R}^{4}$ with
\[
  \star e_{kl} = \varepsilon_{kl}\, e_{k'l'},
\]
that is, it sends the unit area element of each coordinate plane to that of its
orthogonal complement, up to sign:
\[
  \star e_{12} = e_{34}, \quad \star e_{34} = e_{12}, \quad
  \star e_{14} = e_{23}, \quad \star e_{23} = e_{14}, \quad
  \star e_{13} = -e_{24}, \quad \star e_{24} = -e_{13}.
\]
It permutes the orthonormal basis $e_{kl}$ up to signs, so it is an isometry.
Its use here is that, by Laplace expansion, it turns a $4 \times 4$
determinant into an inner product of bivectors.

\subsection{The Laplace expansion}

\begin{theorem}[{The Laplace expansion theorem, \cite[Section~0.8.9]{horn-johnson-2013}}]
  \label{thm:laplace}
  Let $A \in M_{n}(F)$, let $k \in \{1, \ldots, n\}$ be given, and let
  $\beta \subseteq \{1, \ldots, n\}$ be any given index set of cardinality $k$.
  Then
  \[
    \det A = \sum_{\alpha} (-1)^{p(\alpha,\beta)} \det A[\alpha, \beta]
      \det A[\alpha^{c}, \beta^{c}]
    = \sum_{\alpha} (-1)^{p(\alpha,\beta)} \det A[\beta, \alpha]
      \det A[\beta^{c}, \alpha^{c}]
  \]
  in which the sums are over all index sets $\alpha \subseteq \{1, \ldots, n\}$
  of cardinality $k$, and
  $p(\alpha, \beta) = \sum_{i \in \alpha} i + \sum_{j \in \beta} j$.
\end{theorem}

Here $M_{n}(F)$ denotes the $n \times n$ matrices over a field $F$,
$A[\alpha, \beta]$ is the submatrix of $A$ with rows indexed by $\alpha$ and
columns indexed by $\beta$, and $\alpha^{c}$, $\beta^{c}$ are the complements
of $\alpha$, $\beta$ \cite[Section~0.7.1]{horn-johnson-2013}. We use the first
expansion of Theorem~\ref{thm:laplace} for the $4 \times 4$ matrix $A$ with columns $v, v', w, w'$, index
sets of cardinality $2$, and $\beta = \{1, 2\}$, the columns of $v, v'$. For
$\alpha = \{k, l\}$, $k < l$, with complement $\{k', l'\}$, $k' < l'$, the
minor $\det A[\alpha, \beta]$ is $v_{k}v'_{l} - v_{l}v'_{k} =
(v \wedge v')_{kl}$, the complementary minor $\det A[\alpha^{c}, \beta^{c}]$
is $(w \wedge w')_{k'l'}$, and $(-1)^{p(\alpha, \beta)} = (-1)^{k + l + 1 + 2}$,
which is the sign $\varepsilon_{kl}$ of the permutation $(k, l, k', l')$.
Hence
\[
  \det(v, v', w, w')
  = \sum_{k<l} \varepsilon_{kl}\,(v \wedge v')_{kl}\,(w \wedge w')_{k'l'}.
\]
This sum is an inner product, as the following remark shows.

\begin{remark}
  \label{rem:starcoord}
  For $Y \in \Lambda^{2}\mathbb{R}^{4}$, write $Y = \sum_{k<l} Y_{kl}\, e_{kl}$.
  Linearity of $\star$ gives
  \[
    \star Y = \sum_{k<l} Y_{kl}\, \star e_{kl}
    = \sum_{k<l} \varepsilon_{kl}\, Y_{kl}\, e_{k'l'}.
  \]
  The pair $\{k, l\}$ is the complement of its complement, so the only term of
  this sum that contributes to $e_{kl}$ is the one indexed by $\{k', l'\}$, and
  $(\star Y)_{kl} = \varepsilon_{k'l'} Y_{k'l'}$. Moreover $\varepsilon_{k'l'} =
  \varepsilon_{kl}$, since $(k', l', k, l)$ arises from $(k, l, k', l')$ by
  moving each of $k', l'$ past both $k$ and $l$, that is, by four
  transpositions. Hence $(\star Y)_{kl} = \varepsilon_{kl} Y_{k'l'}$, and with
  $\langle X, Y\rangle = \sum_{k<l} X_{kl} Y_{kl}$,
  \[
    \langle v \wedge v', \star(w \wedge w')\rangle
    = \sum_{k<l} (v \wedge v')_{kl}\, \bigl(\star(w \wedge w')\bigr)_{kl}
    = \sum_{k<l} \varepsilon_{kl}\,(v \wedge v')_{kl}\,(w \wedge w')_{k'l'}.
  \]
\end{remark}

By Remark~\ref{rem:starcoord}, the expansion above becomes
\begin{equation}
  \label{eq:detstar}
  \det(v, v', w, w') = \langle v \wedge v', \star(w \wedge w')\rangle.
\end{equation}
For example, $\det(e_{1}, e_{2}, e_{3}, e_{4}) = \langle e_{12}, \star
e_{34}\rangle = \langle e_{12}, e_{12}\rangle = 1$. So, up to the areas
$|v \wedge v'|$ and $|w \wedge w'|$, the determinant compares the plane of
$v, v'$ with the orthogonal complement of that of $w, w'$.

\begin{remark}
  The description of $\star$ on coordinate planes holds for every plane: for
  a $2$-plane $\Pi \subseteq \mathbb{R}^{4}$ through the origin,
  \begin{equation}
    \label{eq:starplane}
    \star\,\omega_{\Pi} = \pm\,\omega_{\Pi^{\perp}}.
  \end{equation}
  Indeed, extend an orthonormal basis $u_{1}, u_{2}$ of $\Pi$ by one $u_{3},
  u_{4}$ of $\Pi^{\perp}$. By \eqref{eq:detstar}, $\langle u_{i} \wedge u_{j},
  \star(u_{1} \wedge u_{2})\rangle = \det(u_{i}, u_{j}, u_{1}, u_{2})$, which
  vanishes unless $\{i, j\} = \{3, 4\}$ and is $\pm 1$ then. The six
  $u_{i} \wedge u_{j}$, $i < j$, are orthonormal, hence a basis of
  $\Lambda^{2}\mathbb{R}^{4}$, and these inner products are the coordinates
  of $\star\omega_{\Pi}$ in it; so $\star\omega_{\Pi} = \pm\, u_{3} \wedge
  u_{4} = \pm\,\omega_{\Pi^{\perp}}$.
\end{remark}

\section{Lemmas on the join of two cycles}
\label{sec:lemmas}

For $m \geq 3$ and $x = (x_{1}, \ldots, x_{m}) \in (\mathbb{R}^{N})^{m}$, with
indices modulo $m$, let
\begin{equation}
  \label{eq:areavector}
  \mathcal{A}(x) = \sum_{j=1}^{m} x_{j} \wedge x_{j+1}.
\end{equation}
Thus $\frac12\mathcal{A}(x)$ is the \emph{area vector} of the closed polygon
$x_{1} \cdots x_{m}$; for $N = 2$ its single coordinate is the signed area given
by the shoelace formula.

\begin{definition}
  \label{def:join}
  For $m, l \geq 3$ let $\mathcal{J}_{m,l} = (m + l, \mathcal{F}_{m,l})$ with
  \[
    \mathcal{F}_{m,l} = \bigl\{(i,\ i+1,\ m+j,\ m+j+1) : i \in \mathbb{Z}_{m},\
      j \in \mathbb{Z}_{l}\bigr\},
  \]
  where the first two entries are read modulo $m$ in $\{1, \ldots, m\}$ and the
  last two modulo $l$ in $\{m+1, \ldots, m+l\}$. For
  $p \in (\mathbb{R}^{4})^{m+l}$ write $a_{i} = p_{i}$, $b_{j} = p_{m+j}$, and
  \[
    A = \mathcal{A}(a_{1}, \ldots, a_{m}), \qquad
    B = \mathcal{A}(b_{1}, \ldots, b_{l}).
  \]
\end{definition}

Each facet of $\mathcal{J}_{m,l}$ consists of an edge of the cycle $(1, \ldots,
m)$ and an edge of the cycle $(m+1, \ldots, m+l)$, so as a simplicial complex
$\mathcal{J}_{m,l}$ is the join of these two cycles. By
Lemma~\ref{lem:freesum} it is the boundary complex of the direct sum of two
regular polygons; the polytope itself is this direct sum, not a join of
polytopes.

By \eqref{eq:signed-volume}
\begin{equation}
  \label{eq:Djoin}
  D(p) = \sum_{i \in \mathbb{Z}_{m}}\sum_{j \in \mathbb{Z}_{l}}
    \det(a_{i}, a_{i+1}, b_{j}, b_{j+1}).
\end{equation}

\begin{lemma}
  \label{lem:freesum}
  Let $m, l \geq 3$, let $u_{1}, u_{2}, u_{3}, u_{4}$ be an orthonormal basis of
  $\mathbb{R}^{4}$, and for $i \in \mathbb{Z}_{m}$, $j \in \mathbb{Z}_{l}$ let
  \[
    a_{i} = \cos\tfrac{2\pi i}{m}\,u_{1} + \sin\tfrac{2\pi i}{m}\,u_{2}, \qquad
    b_{j} = \cos\tfrac{2\pi j}{l}\,u_{3} + \sin\tfrac{2\pi j}{l}\,u_{4}.
  \]
  Then $P = \conv\{a_{i}, b_{j}\} \in \mathcal{P}_{4}(m + l)$, and the point
  $p = (a_{1}, \ldots, a_{m}, b_{1}, \ldots, b_{l})$ realises
  $\mathcal{J}_{m,l}$.
\end{lemma}

\begin{proof}
  In the coordinates of the basis $u_{1}, \ldots, u_{4}$, $P$ is the direct sum
  of two regular polygons with $m$ and $l$ vertices and edges, inscribed in unit
  circles centred at the origin. By the cited statement on direct sums, $P$ is
  $4$-dimensional, all $m + l$ points $a_{i}, b_{j}$ are vertices, and $P$ has
  $ml$ facets. All the points are unit vectors, so $P \in \mathcal{P}_{4}(m +
  l)$.

  It remains to show that each of the $ml$ distinct simplices
  $\conv\{a_{i}, a_{i+1}, b_{j}, b_{j+1}\}$ is a facet. Fix $i \in \mathbb{Z}_{m}$
  and $j \in \mathbb{Z}_{l}$. Let
  \[
    \nu_{i} = \cos\tfrac{(2i+1)\pi}{m}\,u_{1} + \sin\tfrac{(2i+1)\pi}{m}\,u_{2},
    \qquad h_{m} = \cos\tfrac{\pi}{m},
  \]
  so $\nu_{i}$ is the unit vector towards the midpoint of the edge from
  $a_{i}$ to $a_{i+1}$, and $h_{m}$ is the distance of that edge from the
  origin. For every $k$,
  \[
    \nu_{i} \cdot a_{k} = \cos\tfrac{(2k-2i-1)\pi}{m} \leq h_{m},
  \]
  with equality if and only if $2k - 2i - 1 \equiv \pm 1 \pmod{2m}$, that is,
  $k \in \{i, i+1\}$; and $\nu_{i} \cdot b_{k} = 0$, since $\nu_{i}$ lies in
  the plane of $u_{1}, u_{2}$ and $b_{k}$ in that of $u_{3}, u_{4}$. Define
  $\mu_{j}$ and $h_{l} = \cos\frac{\pi}{l}$ in the same way for the edge from
  $b_{j}$ to $b_{j+1}$, and let
  \[
    \varphi(x) = \frac{\nu_{i} \cdot x}{h_{m}} + \frac{\mu_{j} \cdot x}{h_{l}}.
  \]
  At a vertex $a_{k}$ the second term vanishes and the first is at most $1$,
  with equality exactly for $k \in \{i, i+1\}$; at a vertex $b_{k}$ the first
  term vanishes and the second is at most $1$, with equality exactly for
  $k \in \{j, j+1\}$. A linear function attains its maximum over $P$ at
  vertices, so $\varphi \leq 1$ on $P$, and $\{\varphi = 1\}$ is a supporting
  hyperplane with
  \[
    P \cap \{\varphi = 1\} = \conv\{a_{i}, a_{i+1}, b_{j}, b_{j+1}\}.
  \]
  This face is $3$-dimensional: $a_{i}, a_{i+1}$ span the plane of
  $u_{1}, u_{2}$ and $b_{j}, b_{j+1}$ that of $u_{3}, u_{4}$, so the four
  points are linearly, hence affinely, independent. Thus each of the $ml$
  simplices is a facet, and since $P$ has exactly $ml$ facets, these are all
  of them; so $p$ realises $\mathcal{J}_{m,l}$.
\end{proof}

\begin{lemma}
  \label{lem:coherent}
  The orientation of $\mathcal{J}_{m,l}$ is coherent.
\end{lemma}

\begin{proof}
  Recall that deleting entry $a$ of a facet induces the orientation of the
  remaining triple if $a$ is even and the opposite one if $a$ is odd. A ridge
  consists of three of the four vertices of a facet
  $(a_{i}, a_{i+1}, b_{j}, b_{j+1})$, so it is of one of two kinds.
  A facet containing the ridge $\{a_{i}, b_{j}, b_{j+1}\}$ contains the edge
  $\{b_{j}, b_{j+1}\}$ of the second cycle and an edge of the first cycle
  through $a_{i}$, that is, $\{a_{i-1}, a_{i}\}$ or $\{a_{i}, a_{i+1}\}$. So
  this ridge lies in exactly the facets $(a_{i-1}, a_{i}, b_{j}, b_{j+1})$ and
  $(a_{i}, a_{i+1}, b_{j}, b_{j+1})$, which are different since $m \geq 3$.
  Deleting entry~$1$ of the first, an odd position, gives the opposite of
  $(a_{i}, b_{j}, b_{j+1})$; deleting entry~$2$ of the second, an even
  position, gives $(a_{i}, b_{j}, b_{j+1})$ itself.
  In the same way the ridge $\{a_{i}, a_{i+1}, b_{j}\}$ lies in exactly
  $(a_{i}, a_{i+1}, b_{j-1}, b_{j})$ and $(a_{i}, a_{i+1}, b_{j}, b_{j+1})$,
  which are different since $l \geq 3$. Deleting entry~$3$ of the first gives
  the opposite of $(a_{i}, a_{i+1}, b_{j})$, and deleting entry~$4$ of the
  second gives $(a_{i}, a_{i+1}, b_{j})$ itself.
  In both cases the ridge lies in exactly two facets, which induce opposite
  orientations on it.
\end{proof}

\begin{lemma}
  \label{lem:factor}
  For every $p \in (\mathbb{R}^{4})^{m+l}$, with $D(p)$ as in \eqref{eq:Djoin}
  and $A$, $B$ as in Definition~\ref{def:join},
  $D(p) = \langle A, \star B\rangle$. Hence $|D(p)| \leq |A|\,|B|$, and if
  $A, B \neq 0$, equality holds if and only if $\star B = \lambda A$ for some
  $\lambda \in \mathbb{R}$.
\end{lemma}

\begin{proof}
  Applying \eqref{eq:detstar} with $v = a_{i}$, $v' = a_{i+1}$, $w = b_{j}$
  and $w' = b_{j+1}$, each term of \eqref{eq:Djoin} is
  \[
    \det(a_{i}, a_{i+1}, b_{j}, b_{j+1})
    = \langle a_{i} \wedge a_{i+1}, \star(b_{j} \wedge b_{j+1})\rangle.
  \]
  The inner product is bilinear and $\star$ is linear, so summing over $i$
  and $j$ gives
  \[
    D(p) = \Bigl\langle \sum_{i \in \mathbb{Z}_{m}} a_{i} \wedge a_{i+1},\
      \star \sum_{j \in \mathbb{Z}_{l}} b_{j} \wedge b_{j+1} \Bigr\rangle
    = \langle A, \star B\rangle.
  \]
  By the Cauchy--Schwarz inequality, $|\langle A, \star B\rangle| \leq
  |A|\,|\star B|$, and $|\star B| = |B|$ since $\star$ is an isometry; hence
  $|D(p)| \leq |A|\,|B|$. If $A, B \neq 0$, then also $\star B \neq 0$, and
  equality in the Cauchy--Schwarz inequality for two nonzero vectors holds
  if and only if they are linearly dependent, that is, $\star B = \lambda A$
  for some $\lambda \in \mathbb{R}$.
\end{proof}

\section{An inequality for closed polygons}
\label{sec:polygon}

By Lemma~\ref{lem:factor}, $|D(p)| \leq |A|\,|B|$. So it suffices to bound
$|\mathcal{A}(x)|$, defined in \eqref{eq:areavector}, for a closed polygon
with unit vertices $x_{1}, \ldots, x_{m} \in \mathbb{R}^{N}$, and to determine
when the bound is attained. Such a polygon need be neither planar nor convex.
The map $\mathcal{A}$ is quadratic in $x$ and unchanged by a cyclic shift of
the indices, so we expand $x$ in the real discrete Fourier basis, as
Schoenberg did for plane polygons \cite{schoenberg-1950}.

\begin{notation}
  Let $m \geq 3$ and $\theta = 2\pi/m$, and take indices modulo $m$. Let
  $\mathcal{B} \subset \mathbb{R}^{m}$ consist of the vectors $f_{k}$,
  $g_{k}$ for $1 \leq k < m/2$, with entries
  \[
    f_{k,j} = \sqrt{2/m}\,\cos(jk\theta), \qquad
    g_{k,j} = \sqrt{2/m}\,\sin(jk\theta),
  \]
  of $f_{0}$, with all entries $1/\sqrt{m}$, and, if $m$ is even, of
  $f_{m/2}$, with entries $(-1)^{j}/\sqrt{m}$. For $x = (x_{1}, \ldots,
  x_{m}) \in (\mathbb{R}^{N})^{m}$ let
  \[
    y_{k} = \sum_{j=1}^{m} f_{k,j}\, x_{j}, \qquad
    z_{k} = \sum_{j=1}^{m} g_{k,j}\, x_{j} \qquad \in \mathbb{R}^{N}
  \]
  for every $f_{k}$, respectively $g_{k}$, in $\mathcal{B}$.
\end{notation}

So $f_{k}$ and $g_{k}$ are the cosine and sine of frequency $k$ sampled at
the $m$ vertices, and $f_{0}$ and $f_{m/2}$ are the cosines of frequency $0$
and $m/2$, where the sine vanishes; thus $\mathcal{B}$ has $m$ elements. The
vectors $y_{k}$, $z_{k}$ are the real discrete Fourier coefficients of $x$.
The following lemma says that $x$ is recovered from them and that its
squared length splits over them. In complex form, these are the finite
Fourier inversion and Parseval formulas
\cite[Chapter~7, Section~1]{stein-shakarchi-2003},
\cite[Section~1]{schoenberg-1950}.

\begin{lemma}
  \label{lem:dftbasis}
  The set $\mathcal{B}$ is an orthonormal basis of $\mathbb{R}^{m}$. Hence,
  for every $x \in (\mathbb{R}^{N})^{m}$ and $j = 1, \ldots, m$,
  \begin{equation}
    \label{eq:expansion}
    x_{j} = \tfrac{1}{\sqrt{m}}\, y_{0}
      + \sum_{1 \leq k < m/2} x^{(k)}_{j}
      \;\Bigl[ + \tfrac{(-1)^{j}}{\sqrt{m}}\, y_{m/2} \Bigr],
    \qquad
    x^{(k)}_{j} = \sqrt{2/m}\,\bigl(\cos(jk\theta)\, y_{k}
      + \sin(jk\theta)\, z_{k}\bigr),
  \end{equation}
  and
  \begin{equation}
    \label{eq:parseval}
    \sum_{j=1}^{m} |x_{j}|^{2} = |y_{0}|^{2}
      + \sum_{1 \leq k < m/2} \bigl(|y_{k}|^{2} + |z_{k}|^{2}\bigr)
      \;\bigl[ + |y_{m/2}|^{2} \bigr],
  \end{equation}
  where the bracketed terms appear only if $m$ is even.
\end{lemma}

\begin{proof}
  Steps~1--4 show that $\mathcal{B}$ is orthonormal; since it has $m$
  elements, it is then an orthonormal basis of $\mathbb{R}^{m}$. Step~5
  deduces \eqref{eq:expansion} and \eqref{eq:parseval}.

  \emph{Step 1: sums over the vertices.} Steps~2--4 rely on sums of
  $\cos(jn\theta)$ and $\sin(jn\theta)$ over the $m$ vertices, so we
  evaluate these first. For every integer $n$,
  \begin{equation}
    \label{eq:rootsum}
    \sum_{j=1}^{m} \cos(jn\theta) =
      \begin{cases} m, & m \mid n,\\ 0, & m \nmid n, \end{cases}
    \qquad
    \sum_{j=1}^{m} \sin(jn\theta) = 0.
  \end{equation}
  Geometrically, the points $(\cos(jn\theta), \sin(jn\theta))$,
  $j = 1, \ldots, m$, are vertices of a regular polygon centred at the
  origin, possibly traversed several times, and \eqref{eq:rootsum} says
  that they add up to the origin unless they all coincide.

  If $m \mid n$, then $jn\theta$ is a multiple of $2\pi$ for every $j$, so
  every cosine equals $1$ and every sine equals $0$. Otherwise we combine
  the two sums into one complex sum. Let $\zeta = e^{in\theta}$. By
  Euler's formula $e^{i\alpha} = \cos\alpha + i\sin\alpha$ with
  $\alpha = jn\theta$, we have $\zeta^{j} = \cos(jn\theta) + i\sin(jn\theta)$,
  and summing over $j$,
  \[
    \sum_{j=1}^{m} \zeta^{j}
      = \sum_{j=1}^{m} \cos(jn\theta) + i \sum_{j=1}^{m} \sin(jn\theta).
  \]
  So it suffices to show that the left-hand side vanishes. It is a
  geometric sum, and $\zeta \neq 1$ because $n\theta = 2\pi n/m$ is not a
  multiple of $2\pi$, so the geometric sum formula applies. Since
  $m\theta = 2\pi$, we have $\zeta^{m} = e^{imn\theta} = e^{2\pi i n} = 1$,
  and therefore
  \[
    \sum_{j=1}^{m} \zeta^{j} = \frac{\zeta(\zeta^{m} - 1)}{\zeta - 1} = 0.
  \]

  \emph{Step 2: inner products reduce to these sums.} To check
  orthonormality we need the inner product $u \cdot v = \sum_{j} u_{j}
  v_{j}$ for all $u, v \in \mathcal{B}$. We first bring all vectors of
  $\mathcal{B}$ into a common form, and then write their products as the
  sums of Step~1.

  Call $j \mapsto \cos(jk\theta)$ the \emph{cosine} and
  $j \mapsto \sin(jk\theta)$ the \emph{sine} of frequency $k$. Every vector
  of $\mathcal{B}$ is a positive multiple of one of these with
  $0 \leq k \leq m/2$:
  \begin{itemize}
    \item for $1 \leq k < m/2$, $f_{k}$ and $g_{k}$ are $\sqrt{2/m}$ times
      the cosine and the sine of frequency $k$;
    \item $f_{0}$ is $1/\sqrt{m}$ times the cosine of frequency $0$, since
      $\cos(0) = 1$;
    \item $f_{m/2}$ is $1/\sqrt{m}$ times the cosine of frequency $m/2$,
      since $(m/2)\theta = \pi$ and $\cos(j\pi) = (-1)^{j}$.
  \end{itemize}
  From this list, $\mathcal{B}$ contains at most one cosine and at most one
  sine of each frequency. So, up to a positive factor, the inner product of
  two vectors of $\mathcal{B}$ with frequencies $k$ and $l$ is one of the
  three sums on the left of \eqref{eq:products} below; a cosine of
  frequency $k$ times a sine of frequency $l$ is the third sum with $k$ and
  $l$ swapped.

  To evaluate these sums, apply the product-to-sum formulas
  \[
    \begin{gathered}
      2\cos a\cos b = \cos(a-b) + \cos(a+b), \qquad
      2\sin a\sin b = \cos(a-b) - \cos(a+b), \\
      2\sin a\cos b = \sin(a+b) + \sin(a-b)
    \end{gathered}
  \]
  with $a = jk\theta$ and $b = jl\theta$, so that $a \pm b = j(k \pm l)\theta$,
  and sum over $j$. For all $0 \leq k, l \leq m/2$ this gives
  \begin{equation}
    \label{eq:products}
    \begin{aligned}
      \sum_{j} \cos(jk\theta)\cos(jl\theta)
        &= \tfrac12 \sum_{j} \cos(j(k-l)\theta)
          + \tfrac12 \sum_{j} \cos(j(k+l)\theta), \\
      \sum_{j} \sin(jk\theta)\sin(jl\theta)
        &= \tfrac12 \sum_{j} \cos(j(k-l)\theta)
          - \tfrac12 \sum_{j} \cos(j(k+l)\theta), \\
      \sum_{j} \sin(jk\theta)\cos(jl\theta)
        &= \tfrac12 \sum_{j} \sin(j(k+l)\theta)
          + \tfrac12 \sum_{j} \sin(j(k-l)\theta) = 0.
    \end{aligned}
  \end{equation}
  The right-hand sides are sums of the form \eqref{eq:rootsum} with
  $n = k - l$ and $n = k + l$. The last line vanishes because, by
  \eqref{eq:rootsum}, every sum of sines over the vertices is $0$.

  \emph{Step 3: distinct vectors are orthogonal.} Let $u \neq v$ be vectors
  of $\mathcal{B}$ with frequencies $k$ and $l$. We show $u \cdot v = 0$ by
  considering the two possible types of the pair.

  If one of $u, v$ is a cosine and the other a sine, then by Step~2
  $u \cdot v$ is a positive multiple of the third sum in
  \eqref{eq:products}, which is $0$.

  Otherwise $u$ and $v$ are both cosines or both sines, and $u \cdot v$ is
  a positive multiple of the first or the second sum in
  \eqref{eq:products}. Their right-hand sides consist of cosine sums with
  $n = k - l$ and $n = k + l$, and by \eqref{eq:rootsum} such a sum
  vanishes as soon as $m \nmid n$. So it suffices to show that $m$ divides
  neither $k - l$ nor $k + l$:
  \begin{itemize}
    \item $k \neq l$: since $\mathcal{B}$ has at most one cosine and at most
      one sine of each frequency, two distinct cosines, or two distinct
      sines, have different frequencies;
    \item $0 < |k - l| < m$: the difference is nonzero since $k \neq l$,
      and at most $m/2$ since $0 \leq k, l \leq m/2$;
    \item $0 < k + l < m$: the sum is positive since $k$ and $l$ are not
      both $0$, and less than $m$ since both are at most $m/2$ and they are
      not both equal to $m/2$.
  \end{itemize}
  A nonzero integer strictly between $-m$ and $m$ is not divisible by $m$,
  so both cosine sums vanish, and $u \cdot v = 0$.

  \emph{Step 4: unit length.} It remains to show $|v|^{2} = v \cdot v = 1$
  for every $v \in \mathcal{B}$. We treat $f_{k}, g_{k}$ with
  $1 \leq k < m/2$ using \eqref{eq:products}, and $f_{0}, f_{m/2}$
  directly.

  Let $1 \leq k < m/2$. Since $f_{k,j} = \sqrt{2/m}\,\cos(jk\theta)$ and
  $g_{k,j} = \sqrt{2/m}\,\sin(jk\theta)$,
  \[
    |f_{k}|^{2} = \frac{2}{m} \sum_{j} \cos^{2}(jk\theta),
    \qquad
    |g_{k}|^{2} = \frac{2}{m} \sum_{j} \sin^{2}(jk\theta),
  \]
  so it suffices to show that both sums equal $m/2$. These are the first
  two lines of \eqref{eq:products} with $l = k$:
  \[
    \sum_{j} \cos^{2}(jk\theta)
      = \tfrac12 \sum_{j} \cos(0) + \tfrac12 \sum_{j} \cos(2jk\theta),
    \qquad
    \sum_{j} \sin^{2}(jk\theta)
      = \tfrac12 \sum_{j} \cos(0) - \tfrac12 \sum_{j} \cos(2jk\theta).
  \]
  The first sum on each right-hand side is $\sum_{j} 1 = m$. The second is
  the sum \eqref{eq:rootsum} with $n = 2k$, and it vanishes because
  $0 < 2k < m$, so $m \nmid 2k$. Hence both sums equal $\tfrac12 m$, and
  $|f_{k}|^{2} = |g_{k}|^{2} = \frac{2}{m} \cdot \frac{m}{2} = 1$.

  Finally, each of the $m$ entries of $f_{0}$ and of $f_{m/2}$ is
  $\pm 1/\sqrt{m}$, so each squared entry is $1/m$ and
  $|f_{0}|^{2} = |f_{m/2}|^{2} = m \cdot \frac1m = 1$.

  \emph{Step 5: expansion and Parseval.} We now deduce
  \eqref{eq:expansion} and \eqref{eq:parseval}. By Steps~1--4,
  $\mathcal{B}$ is an orthonormal basis of $\mathbb{R}^{m}$, so every
  $\xi \in \mathbb{R}^{m}$ is the sum of its projections onto the basis
  vectors, and its squared length is the sum of the squared coefficients:
  \begin{equation}
    \label{eq:onb}
    \xi = \sum_{v \in \mathcal{B}} (\xi \cdot v)\, v,
    \qquad
    |\xi|^{2} = \sum_{v \in \mathcal{B}} (\xi \cdot v)^{2}.
  \end{equation}
  The vertices $x_{j}$, however, lie in $\mathbb{R}^{N}$, so we apply
  \eqref{eq:onb} to one coordinate at a time.

  Fix $c \in \{1, \ldots, N\}$ and let $\xi^{(c)} \in \mathbb{R}^{m}$ be the
  vector of $c$-th coordinates of the vertices, $\xi^{(c)}_{j} =
  (x_{j})_{c}$. Taking the $c$-th coordinate in the definition of $y_{k}$
  and $z_{k}$ gives
  \[
    (y_{k})_{c} = \sum_{j} f_{k,j}\, (x_{j})_{c} = \xi^{(c)} \cdot f_{k},
    \qquad
    (z_{k})_{c} = \xi^{(c)} \cdot g_{k}.
  \]
  So for $\xi = \xi^{(c)}$, the coefficients $\xi \cdot v$ in
  \eqref{eq:onb} are the $c$-th coordinates of the vectors $y_{k}$ and
  $z_{k}$.

  For \eqref{eq:expansion}, take the $j$-th entry of the first identity in
  \eqref{eq:onb} with $\xi = \xi^{(c)}$. Its left side is $(x_{j})_{c}$,
  and its right side is the sum of the following terms
  $(\xi^{(c)} \cdot v)\, v_{j}$:
  \begin{itemize}
    \item $v = f_{0}$ gives $(y_{0})_{c}/\sqrt{m}$, since
      $f_{0,j} = 1/\sqrt{m}$;
    \item $v = f_{m/2}$, if $m$ is even, gives
      $(-1)^{j} (y_{m/2})_{c}/\sqrt{m}$, since
      $f_{m/2,j} = (-1)^{j}/\sqrt{m}$;
    \item $v = f_{k}$ and $v = g_{k}$, for $1 \leq k < m/2$, together give
      \[
        (y_{k})_{c}\, f_{k,j} + (z_{k})_{c}\, g_{k,j}
          = \sqrt{2/m}\,\bigl(\cos(jk\theta)\, (y_{k})_{c}
            + \sin(jk\theta)\, (z_{k})_{c}\bigr)
          = \bigl(x^{(k)}_{j}\bigr)_{c}.
      \]
  \end{itemize}
  These are the $c$-th coordinates of the terms on the right of
  \eqref{eq:expansion}. As $c$ is arbitrary, \eqref{eq:expansion} holds.

  For \eqref{eq:parseval}, apply the second identity in \eqref{eq:onb} to
  $\xi^{(c)}$ and sum over $c = 1, \ldots, N$. On the left,
  $\sum_{c} |\xi^{(c)}|^{2} = \sum_{c} \sum_{j} (x_{j})_{c}^{2}
  = \sum_{j} |x_{j}|^{2}$. On the right, each $v \in \mathcal{B}$
  contributes $\sum_{c} (\xi^{(c)} \cdot v)^{2}$. For $v = f_{k}$,
  including $f_{0}$ and $f_{m/2}$, this is
  $\sum_{c} (y_{k})_{c}^{2} = |y_{k}|^{2}$; for $v = g_{k}$ it is
  $|z_{k}|^{2}$. This is \eqref{eq:parseval}.
\end{proof}

The next lemma splits $\mathcal{A}(x)$ over the frequencies; for $N = 2$ it
is due to Schoenberg, in complex form \cite[Section~2]{schoenberg-1950}.

\begin{lemma}
  \label{lem:fourier}
  For every $x \in (\mathbb{R}^{N})^{m}$,
  \begin{equation}
    \label{eq:fourier}
    \mathcal{A}(x) = \sum_{1 \leq k < m/2} 2\sin(k\theta)\; y_{k} \wedge z_{k}.
  \end{equation}
\end{lemma}

\begin{proof}
  The idea is to substitute the expansion \eqref{eq:expansion} into
  $\mathcal{A}(x)$. This works best through the symmetric difference
  $x_{j+1} - x_{j-1}$, which turns each cosine into a sine and each sine
  into a cosine (Step~2), so we first rewrite $\mathcal{A}(x)$ in terms of
  it.

  \emph{Step 1: symmetric form.} Substituting $j + 1$ for $j$, which only
  permutes the terms since indices are taken modulo $m$, and using
  $a \wedge b = -b \wedge a$,
  \[
    \sum_{j} x_{j} \wedge x_{j-1} = \sum_{j} x_{j+1} \wedge x_{j}
      = -\sum_{j} x_{j} \wedge x_{j+1} = -\mathcal{A}(x).
  \]
  Subtracting this from $\mathcal{A}(x) = \sum_{j} x_{j} \wedge x_{j+1}$
  gives
  \[
    \mathcal{A}(x) = \frac12 \sum_{j=1}^{m} x_{j} \wedge (x_{j+1} - x_{j-1}).
  \]

  \emph{Step 2: the symmetric difference.} We compute $x_{j+1} - x_{j-1}$
  term by term from \eqref{eq:expansion}.
  \begin{itemize}
    \item The constant term $y_{0}/\sqrt{m}$ does not depend on $j$, so it
      cancels.
    \item For even $m$, the term $(-1)^{j} y_{m/2}/\sqrt{m}$ takes the same
      value at $j + 1$ and $j - 1$, since $(-1)^{j+1} = (-1)^{j-1}$, so it
      cancels too.
    \item For $1 \leq k < m/2$, the addition formulas
      $\cos(\alpha \pm \varphi) = \cos\alpha\cos\varphi \mp
      \sin\alpha\sin\varphi$ and $\sin(\alpha \pm \varphi) =
      \sin\alpha\cos\varphi \pm \cos\alpha\sin\varphi$ with
      $\alpha = jk\theta$ and $\varphi = k\theta$ give
      \[
        \begin{aligned}
          \cos((j+1)k\theta) - \cos((j-1)k\theta)
            &= -2\sin(k\theta)\,\sin(jk\theta), \\
          \sin((j+1)k\theta) - \sin((j-1)k\theta)
            &= 2\sin(k\theta)\,\cos(jk\theta).
        \end{aligned}
      \]
      Applied to $x^{(k)}_{j} = \sqrt{2/m}\,(\cos(jk\theta)\, y_{k} +
      \sin(jk\theta)\, z_{k})$, this gives
      \[
        x^{(k)}_{j+1} - x^{(k)}_{j-1} = 2\sin(k\theta)\, \sqrt{2/m}\,
          \bigl(\cos(jk\theta)\, z_{k} - \sin(jk\theta)\, y_{k}\bigr).
      \]
  \end{itemize}
  Hence
  \[
    x_{j+1} - x_{j-1} = \sum_{1 \leq k < m/2} 2\sin(k\theta)\,
      \sqrt{2/m}\,\bigl(\cos(jk\theta)\, z_{k} - \sin(jk\theta)\, y_{k}\bigr).
  \]

  \emph{Step 3: summing over the vertices.} Substitute this into Step~1.
  Since $\wedge$ is bilinear, we may pull the scalars out and carry out the
  sum over $j$ first, using $\sum_{j} (a_{j} x_{j}) \wedge z =
  (\sum_{j} a_{j} x_{j}) \wedge z$:
  \[
    \mathcal{A}(x)
      = \sum_{1 \leq k < m/2} \sin(k\theta)\,
        \Bigl( \Bigl(\sqrt{2/m} \sum_{j} \cos(jk\theta)\, x_{j}\Bigr)
          \wedge z_{k}
        - \Bigl(\sqrt{2/m} \sum_{j} \sin(jk\theta)\, x_{j}\Bigr)
          \wedge y_{k} \Bigr).
  \]
  The two bracketed sums are, by definition, $y_{k}$ and $z_{k}$. Using
  $z_{k} \wedge y_{k} = -y_{k} \wedge z_{k}$,
  \[
    \mathcal{A}(x)
      = \sum_{1 \leq k < m/2} \sin(k\theta)\,
        \bigl(y_{k} \wedge z_{k} - z_{k} \wedge y_{k}\bigr)
      = \sum_{1 \leq k < m/2} 2\sin(k\theta)\; y_{k} \wedge z_{k}.
    \qedhere
  \]
\end{proof}

\begin{proposition}
  \label{prop:polygon}
  Let $m \in \{3, 4, 5\}$ and $x \in (\mathbb{R}^{N})^{m}$ with $|x_{j}| = 1$
  for all $j$. Then $|\mathcal{A}(x)| \leq m\sin\frac{2\pi}{m}$, with equality
  if and only if
  \begin{equation}
    \label{eq:regular}
    x_{j} = \cos\frac{2\pi j}{m}\,u + \sin\frac{2\pi j}{m}\,w
      \qquad (j = 1, \ldots, m)
  \end{equation}
  for orthonormal $u, w$, that is, if the $x_{j}$ are the consecutive vertices
  of a regular $m$-gon inscribed in the unit circle of $\spn(u, w)$. In this
  case $\mathcal{A}(x) = m\sin\frac{2\pi}{m}\; u \wedge w$.
\end{proposition}

\begin{proof}
  Recall $\theta = 2\pi/m$, so the claimed bound is $m\sin\theta$. Sums
  over $k$ run over $1 \leq k < m/2$, that is, over $k = 1$ for $m = 3, 4$
  and over $k = 1, 2$ for $m = 5$. Steps~1--4 each prove one inequality,
  starting from Lemma~\ref{lem:fourier}, and record when it is an
  equality; chained together they give the bound. Step~5 shows that
  equality forces a regular polygon, and Step~6 that the regular polygon
  attains the bound.

  \emph{Step 1: triangle inequality.} For $1 \leq k < m/2$ we have
  $0 < k\theta < \pi$, so the weights $\sin(k\theta)$ in \eqref{eq:fourier}
  are positive, and the triangle inequality gives
  \[
    |\mathcal{A}(x)| \leq \sum_{k} 2\sin(k\theta)\, |y_{k} \wedge z_{k}|.
  \]

  \emph{Step 2: Lagrange's identity.} By Lemma~\ref{lem:lagrange} in the
  form stated after it, $2|y_{k} \wedge z_{k}| \leq |y_{k}|^{2} +
  |z_{k}|^{2}$, so
  \[
    \sum_{k} 2\sin(k\theta)\, |y_{k} \wedge z_{k}|
      \leq \sum_{k} \sin(k\theta)\,\bigl(|y_{k}|^{2} + |z_{k}|^{2}\bigr).
  \]
  Since the weights are positive, this is an equality if and only if
  $y_{k} \cdot z_{k} = 0$ and $|y_{k}| = |z_{k}|$ for every $k$.

  \emph{Step 3: the largest weight.} We replace each weight by the largest
  one. For $m = 3, 4$ there is only $k = 1$. For $m = 5$,
  $\sin 2\theta = \sin\frac{4\pi}{5} = \sin\frac{\pi}{5} <
  \sin\frac{2\pi}{5} = \sin\theta$. So $\sin(k\theta) \leq \sin\theta$ for
  every $k$, and
  \[
    \sum_{k} \sin(k\theta)\,\bigl(|y_{k}|^{2} + |z_{k}|^{2}\bigr)
      \leq \sin\theta \sum_{k} \bigl(|y_{k}|^{2} + |z_{k}|^{2}\bigr).
  \]
  The two sides differ by $(\sin\theta - \sin 2\theta)(|y_{2}|^{2} +
  |z_{2}|^{2})$ for $m = 5$, and not at all for $m = 3, 4$. So this is an
  equality if and only if $y_{2} = z_{2} = 0$ when $m = 5$.

  \emph{Step 4: Parseval.} Since the vertices are unit vectors,
  $\sum_{j} |x_{j}|^{2} = m$. By \eqref{eq:parseval}, $m$ is the sum on the
  right of Step~3 plus $|y_{0}|^{2}$, and, for $m = 4$, plus
  $|y_{2}|^{2}$, as then $m/2 = 2$. So
  \[
    \sin\theta \sum_{k} \bigl(|y_{k}|^{2} + |z_{k}|^{2}\bigr)
      \leq m\sin\theta,
  \]
  with equality if and only if $y_{0} = 0$ and, for $m = 4$, $y_{2} = 0$.

  Chaining Steps~1--4 gives the bound:
  \[
    \begin{aligned}
      |\mathcal{A}(x)|
        &\stackrel{\text{Step 1}}{\leq}
          \sum_{k} 2\sin(k\theta)\, |y_{k} \wedge z_{k}|
        \stackrel{\text{Step 2}}{\leq}
          \sum_{k} \sin(k\theta)\,\bigl(|y_{k}|^{2} + |z_{k}|^{2}\bigr) \\
        &\stackrel{\text{Step 3}}{\leq}
          \sin\theta \sum_{k} \bigl(|y_{k}|^{2} + |z_{k}|^{2}\bigr)
        \stackrel{\text{Step 4}}{\leq} m\sin\theta.
    \end{aligned}
  \]

  \emph{Step 5: equality forces a regular polygon.} Suppose
  $|\mathcal{A}(x)| = m\sin\theta$. Then every inequality in the chain is
  an equality. By the equality cases of Steps~3 and~4, every Fourier
  coefficient except $y_{1}$ and $z_{1}$ vanishes: $y_{0}$ in all cases,
  $y_{2}$ for $m = 4$, and $y_{2}, z_{2}$ for $m = 5$. So
  \eqref{eq:expansion} reduces to
  \[
    x_{j} = \sqrt{2/m}\,\bigl(\cos(j\theta)\, y_{1} + \sin(j\theta)\,
      z_{1}\bigr),
  \]
  and \eqref{eq:parseval} to $|y_{1}|^{2} + |z_{1}|^{2} = m$. By the
  equality case of Step~2 for $k = 1$, $y_{1} \cdot z_{1} = 0$ and
  $|y_{1}| = |z_{1}|$, hence $|y_{1}| = |z_{1}| = \sqrt{m/2}$. So
  $u = \sqrt{2/m}\, y_{1}$ and $w = \sqrt{2/m}\, z_{1}$ are orthonormal, and
  the last display is \eqref{eq:regular}.

  \emph{Step 6: the regular polygon attains the bound.} Suppose
  \eqref{eq:regular} holds. Expand $x_{j} \wedge x_{j+1}$ by bilinearity;
  since $u \wedge u = w \wedge w = 0$ and $w \wedge u = -u \wedge w$, only
  two terms survive, and by the subtraction formula for the sine,
  \[
    x_{j} \wedge x_{j+1}
      = \bigl(\cos(j\theta)\sin((j+1)\theta)
        - \sin(j\theta)\cos((j+1)\theta)\bigr)\, u \wedge w
      = \sin\theta\; u \wedge w
  \]
  for every $j$. Summing over the $m$ values of $j$ gives
  $\mathcal{A}(x) = m\sin\theta\; u \wedge w$. Finally, $|u \wedge w| = 1$
  by Lemma~\ref{lem:lagrange}, since $|u| = |w| = 1$ and $u \cdot w = 0$, so
  $|\mathcal{A}(x)| = m\sin\theta$.
\end{proof}

\begin{remark}
  \label{rem:planar}
  For $N = 2$, Proposition~\ref{prop:polygon} is the case $m \leq 5$ of
  \cite[Theorem~5]{schoenberg-1950}, applied to $x$ and to its reverse
  (see also \cite[Section~2.3, Theorem~12]{kazarinoff-1961}), and the proof
  above follows Schoenberg's, with one further step,
  Lemma~\ref{lem:lagrange}.
\end{remark}

\begin{remark}
  \label{rem:projection}
  For $m = 3, 4$ the planar case already gives every $N$: by
  Lemma~\ref{lem:fourier}, $\mathcal{A}(x) = 2\sin\theta\; y_{1} \wedge
  z_{1}$ lies in the plane of $y_{1}, z_{1}$, and projecting $x$ onto that
  plane keeps $\mathcal{A}(x)$ and does not increase $\sum_{j}
  |x_{j}|^{2}$. For $m = 5$ this fails, which is why the case needed by the
  join $\mathcal{J}_{4,5}$ requires the argument above: the unit vectors
  $x_{j} = \frac{1}{\sqrt{2}}(\cos(j\theta)\, e_{1} + \sin(j\theta)\, e_{2}
  + \cos(2j\theta)\, e_{3} + \sin(2j\theta)\, e_{4})$ have
  $\mathcal{A}(x) = \frac52(\sin\theta\; e_{12} + \sin 2\theta\; e_{34})$,
  which is not a multiple of $\omega_{\Pi}$ for any plane $\Pi$, so no
  projection onto a plane keeps it.
\end{remark}

\begin{remark}
  The restriction $m \leq 5$ is needed: for $m \geq 7$, the polygon
  $x = x^{(2)}$ with $y_{2}, z_{2}$ orthogonal of length $\sqrt{m/2}$ has
  unit vertices and $|\mathcal{A}(x)| = m\sin 2\theta > m\sin\theta$. For
  $N = 2$, Schoenberg determines the extremal polygons for every $m$
  \cite[Theorem~5]{schoenberg-1950}.
\end{remark}

\section{The optimum within the join type}
\label{sec:pair}

\begin{theorem}
  \label{thm:pair}
  Let $m, l \in \{3, 4, 5\}$, $p \in (\mathbb{S}^{3})^{m+l}$, and $D(p)$ the
  signed volume \eqref{eq:Djoin} of $\mathcal{J}_{m,l}$, with $a_{i}$, $b_{j}$
  as in Definition~\ref{def:join}. Then
  \begin{equation}
    \label{eq:pairbound}
    |D(p)| \leq m\sin\frac{2\pi}{m}\cdot l\sin\frac{2\pi}{l},
  \end{equation}
  with equality if and only if, for some $2$-plane $\Pi$, the points $a_{1},
  \ldots, a_{m}$ are the consecutive vertices of a regular $m$-gon inscribed in
  the unit circle of $\Pi$ and $b_{1}, \ldots, b_{l}$ those of a regular
  $l$-gon inscribed in the unit circle of $\Pi^{\perp}$.
\end{theorem}

\begin{proof}
  With $A$ and $B$ as in Definition~\ref{def:join},
  Lemma~\ref{lem:factor} gives $|D(p)| \leq |A|\,|B|$, and by
  Proposition~\ref{prop:polygon}, $|A| \leq m\sin\frac{2\pi}{m}$ and $|B| \leq
  l\sin\frac{2\pi}{l}$. This proves \eqref{eq:pairbound}.

  If equality holds, then $A, B \neq 0$, since the right side of
  \eqref{eq:pairbound} is positive, and equality holds in all three
  inequalities. By Proposition~\ref{prop:polygon}, the $a_{i}$ are as in
  \eqref{eq:regular} in a plane $\Pi$, with $A = m\sin\frac{2\pi}{m}\,
  \omega_{\Pi}$, and the $b_{j}$ in a plane $\Pi'$, with $B =
  l\sin\frac{2\pi}{l}\,\omega_{\Pi'}$. By Lemma~\ref{lem:factor}, $\star B =
  \lambda A$, so $\star\omega_{\Pi'}$ is a multiple of $\omega_{\Pi}$; both
  have length $1$, as $\star$ is an isometry, so by \eqref{eq:starplane}
  $\omega_{\Pi'^{\perp}} = \pm\star\omega_{\Pi'} = \pm\omega_{\Pi}$. With an
  orthonormal basis $u_{1}, u_{2}$ of $\Pi$, $\omega_{\Pi} = u_{1}u_{2}^{\mathsf{T}} -
  u_{2}u_{1}^{\mathsf{T}}$ has column space $\Pi$, and likewise for
  $\Pi'^{\perp}$; hence $\Pi = \Pi'^{\perp}$.

  Conversely, let $u_{1}, u_{2}$ and $u_{3}, u_{4}$ be orthonormal bases of
  $\Pi$ and $\Pi^{\perp}$ for which the $a_{i}$ and $b_{j}$ have the form
  \eqref{eq:regular}. By Proposition~\ref{prop:polygon}, $A = m\sin\frac{2\pi}{m}\,u_{1}
  \wedge u_{2}$ and $B = l\sin\frac{2\pi}{l}\,u_{3} \wedge u_{4}$, and by
  Lemma~\ref{lem:factor} and \eqref{eq:starplane}, $|D(p)| = |\langle A, \star
  B\rangle| = |A|\,|B| = m\sin\frac{2\pi}{m}\cdot l\sin\frac{2\pi}{l}$.
\end{proof}

\begin{corollary}
  \label{cor:pair}
  Let $m, l \in \{3, 4, 5\}$. If $P \in \mathcal{P}_{4}(m+l)$ is
  combinatorially equivalent to the polytope of Lemma~\ref{lem:freesum}, then
  \[
    \vol_{4}(P) \leq \frac{1}{24}\, m\sin\frac{2\pi}{m}\; l\sin\frac{2\pi}{l}
      = \frac{1}{6}\area(R_{m})\area(R_{l}),
  \]
  with equality if and only if $P$ is the polytope of Lemma~\ref{lem:freesum}
  for some orthonormal basis $u_{1}, \ldots, u_{4}$, that is, up to an
  orthogonal transformation, the convex hull of $R_{m}$ and $R_{l}$ inscribed
  in $\mathbb{S}^{3}$ in orthogonal $2$-planes through the origin.
\end{corollary}

\begin{proof}
  Label the vertices of $P$ so that $p \in (\mathbb{S}^{3})^{m+l}$ realises
  $\mathcal{J}_{m,l}$. By Lemmas~\ref{lem:coherent} and~\ref{lem:volume},
  $24\vol_{4}(P) = |D(p)|$, and Theorem~\ref{thm:pair} gives the bound. In the
  case of equality, by Theorem~\ref{thm:pair} and \eqref{eq:regular} the
  $a_{i}$ and $b_{j}$ have the form of Lemma~\ref{lem:freesum} for an
  orthonormal basis $u_{1}, u_{2}$ of $\Pi$ and one $u_{3}, u_{4}$ of
  $\Pi^{\perp}$, so $P$ is the polytope of that lemma. Conversely, by
  Lemma~\ref{lem:freesum} the polytope of that lemma lies in
  $\mathcal{P}_{4}(m+l)$ and its vertices, in the order given there, realise
  $\mathcal{J}_{m,l}$; they have the form of the equality case of
  Theorem~\ref{thm:pair}, so by Lemmas~\ref{lem:coherent}
  and~\ref{lem:volume} its volume is
  $\frac{1}{24}\, m\sin\frac{2\pi}{m}\; l\sin\frac{2\pi}{l}$.
\end{proof}

\section{The other combinatorial types}
\label{sec:rivals}

Firsching's data \cite{firsching-2017} give a polytope with rational
vertices for each line of a file with $1142$ lines \artref{data}; we label
them \texttt{firsching-9-0} to \texttt{firsching-9-1141}.

\begin{lemma}
  \label{lem:data}
  The polytopes \texttt{firsching-9-0} to \texttt{firsching-9-1141} represent
  every combinatorial type of simplicial $4$-polytope with nine vertices, each
  exactly once, and the type of \texttt{firsching-9-10} is that of the
  polytope of Lemma~\ref{lem:freesum} with $(m, l) = (4, 5)$.
\end{lemma}

\begin{proof}
  We check in exact arithmetic, as in \cite[Remark~2.14]{medgyes-maxvol-d4}
  for eight vertices, that these $1142$ polytopes are simplicial and pairwise
  not combinatorially equivalent \artref{testnine}. By
  Theorem~\ref{thm:count} they represent every type. The facets of
  \texttt{firsching-9-10} are exactly those of $\mathcal{J}_{4,5}$ for the
  $4$-cycle $(1, 3, 6, 4)$ and the $5$-cycle $(0, 2, 8, 7, 5)$ of its
  zero-based vertex labels \artref{testnine}. By Lemma~\ref{lem:freesum}, the
  facets of the polytope there are those of $\mathcal{J}_{4,5}$ too, so the two
  are combinatorially equivalent.
\end{proof}

\begin{proposition}
  \label{prop:rivals}
  Every simplicial $P \in \mathcal{P}_{4}(9)$ that is not combinatorially
  equivalent to \texttt{firsching-9-10} has $\vol_{4}(P) \leq 39/50$.
\end{proposition}

\begin{proof}
  By Lemma~\ref{lem:data}, $P$ is combinatorially equivalent to
  \texttt{firsching-9-$k$} for some $k \neq 10$. For each of these $1141$
  values of $k$, the stored
  certificate with $\beta = 468/25$ passes the checks of
  \cite[Section~3.3]{medgyes-maxvol-d4}, including that its facets are those
  of \texttt{firsching-9-$k$} \artref{certificates}, \artref{verifier},
  \artref{testnine}. For that $k$, let $\mathcal{C}$ be the coherently
  oriented complex for which the certificate is checked. Since $P$ is
  combinatorially equivalent to \texttt{firsching-9-$k$}, its vertices can be
  listed as $p = (p_{1}, \ldots, p_{9})$ so that $p$ realises $\mathcal{C}$,
  and Proposition~\ref{prop:certificate} gives
  $\vol_{4}(P) \leq \beta/4! = 39/50$.
\end{proof}

\section{The main theorem}
\label{sec:main}

\begin{theorem}
  \label{thm:main}
  \[
    v_{4}(9) = \frac{5}{6}\sin\frac{2\pi}{5}
      = \frac{5}{24}\sqrt{10 + 2\sqrt{5}} = 0.792547\ldots,
  \]
  and every $P \in \mathcal{P}_{4}(9)$ with $\vol_{4}(P) = v_{4}(9)$ is, up to
  an orthogonal transformation, the convex hull of a square and a regular
  pentagon inscribed in $\mathbb{S}^{3}$ in orthogonal $2$-planes through the
  origin.
\end{theorem}

\begin{proof}
  Let $V^{*} = \frac{1}{24}\cdot 4\sin\frac{2\pi}{4}\cdot 5\sin\frac{2\pi}{5} =
  \frac{5}{6}\sin\frac{2\pi}{5}$; since $\sin\frac{2\pi}{5} =
  \frac14\sqrt{10 + 2\sqrt5}$, this is the stated value. By
  Lemma~\ref{lem:freesum} and Corollary~\ref{cor:pair} with $(m, l) = (4, 5)$,
  the square--pentagon pair lies in $\mathcal{P}_{4}(9)$ and has volume
  $V^{*}$, so $v_{4}(9) \geq V^{*}$.

  Let $P \in \mathcal{P}_{4}(9)$ with $\vol_{4}(P) = v_{4}(9)$, which exists by
  Remark~\ref{rem:attained}. By Remark~\ref{rem:propertyz} and
  Lemma~\ref{lem:simplicial}, $P$ is simplicial. Since $39/50 < V^{*}$,
  Proposition~\ref{prop:rivals} shows that $P$ has the combinatorial type of
  \texttt{firsching-9-10}, which by Lemma~\ref{lem:data} is that of the
  polytope of Lemma~\ref{lem:freesum} with $(m, l) = (4, 5)$. Corollary~\ref{cor:pair} gives
  $v_{4}(9) = \vol_{4}(P) \leq V^{*}$, so $v_{4}(9) = V^{*}$, and its equality
  case identifies $P$.
\end{proof}

\section*{AI assistance}

An AI assistant was used during this work, including for mathematical
discussion, drafting and revising the text, and checking the cited statements
against their sources. The proofs do not depend on the assistant: they rest on
the cited results, on the arguments given here, and on the certificates, which
the exact verifier checks.

\section*{Code and data}
\label{sec:artifacts}

{\footnotesize
Everything reported here is available at \repohref, in the version tagged
\href{\repobase/tree/\repocommit}{\texttt{\repocommit}}. The entries below are cited
as [R$k$]; they are artifacts rather than literature and are numbered separately
from the bibliography.

\begin{artifactlist}
  \artifact{data}{data/polytopalsimplicialspheres9_4.txt}{Firsching's
    realisations of the $1142$ simplicial types with nine vertices, fetched by
    \artpath{scripts/setup\_data.py}}
  \artifact{firsching}{src/inscribed/enumeration/firsching.py}{reads the data
    file, keeps the inscribed realisations and builds their facets}
  \artifact{certificates}{results/certificates/d4-n09}{the exact certificates
    with $\beta = 468/25$, produced by \artpath{scripts/eliminate.py} and, for
    \texttt{firsching-9-323} and \texttt{firsching-9-324}, by
    \artpath{scripts/eliminate\_non\_inscribable.py}}
  \artifact{verifier}{src/inscribed/certify/verify.py}{the exact verifier of
    \cite[Section~3.3]{medgyes-maxvol-d4}, with
    \artpath{src/inscribed/combinatorics/complex.py} checking coherence}
  \artifact{testnine}{tests/papers/maxvol\_d4\_n9/test\_nine.py}{tests of the
    claims about Firsching's data and the certificates; the checks over all
    $1142$ types run with \texttt{poetry run pytest -m slow}}
\end{artifactlist}
}

\bibliographystyle{amsplain}
\bibliography{references}

\end{document}